\documentclass[11pt, a4paper]{article}

\usepackage{amsmath,amssymb, amsthm}
\usepackage[T2A]{fontenc}
\usepackage[utf8]{inputenc}
\usepackage{mathrsfs}
\usepackage{fullpage}
\usepackage[affil-it]{authblk}
\usepackage{tikz}

\usepackage{verbatim}

\usepackage{hyperref}

\usepackage{xcolor}
\hypersetup{
  colorlinks   = true, 
  urlcolor     = {blue!90!black}, 
  linkcolor    = {blue!90!black}, 
  citecolor   = {red!90!black} 
}

\newcommand{\R}{\mathbb{R}}
\newcommand{\N}{\mathbb{N}}
\newcommand{\Z}{\mathbb{Z}}

\newcommand{\B}{\mathbf{B}}
\newcommand{\X}{\mathcal{X}}
\newcommand{\0}{\mathbf{0}}
\renewcommand{\d}{\mathrm{d}} 

\newcommand{\x}{\mathbf{x}}

\newcommand{\E}{\mathbf{E}}
\newcommand{\e}{\mathbf{e}}

\newcommand\blfootnote[1]{
  \begingroup
  \renewcommand\thefootnote{}\footnote{#1}%
  \addtocounter{footnote}{-1}%
  \endgroup
}

\theoremstyle{plain}
\newtheorem{thm}{Theorem}[section]

\newtheorem{lem}[thm]{Lemma}

\theoremstyle{definition}

\newtheorem{rmk}[thm]{Remark}

\author{Viktor Bezborodov \thanks{Email: \texttt{vbezborodov@math.uni-bielefeld.de}} }

\affil{\emph{Bielefeld University, Faculty of Mathematics}}

\title{An exponential estimate for the extinction time 
of the branching random walk on a cube}

\begin{document}
 
%
%

\maketitle

\begin{abstract}

 We prove the exponential estimate 
 \begin{equation*} 
 P \{ s < \tau < \infty \} \leq C e^{-q s},   \quad s \geq 0,
\end{equation*}
where $C, q >0$ are constants and $ \tau $ is 
the extinction time of the supercritical branching random walk (BRW) on a cube.
We cover both the discrete-space and continuous-space BRWs.

\end{abstract}


\section{Introduction}

\blfootnote{Keywords: \emph{branching random walk,
exponential estimate,
oriented percolation}}

In this short paper 
we prove an exponential estimate 
for the extinction time of 
a 
branching random walk on a  cube. 
We treat both the discrete-space and
continuous-space models. 
Time is continuous in both models.
A detailed description of 
them can be found in Section 2.

More specifically, we prove the
 exponential estimate
\begin{equation} \label{exp est}
 P \{ s < \tau < \infty \} \leq C e^{-q s},   \quad s \geq 0,
\end{equation}
where $C, q > 0$ are some constants and $ \tau $ is 
the extinction time.
For supercritical spatial random structures, first
estimates of this type
have probably been obtained 
for the oriented percolation process in two dimensions,
see Durrett \cite{Dur84};
for the supercritical contact process,
see e.g. Theorem 2.30 in Liggett \cite{Lig99}.

This work relies on results of
Mountford and Schinazi \cite{MS05} and
Bertacchi and Zucca \cite{BZ09} (see also \cite{BZ15}), who  proved 
in discrete-space settings
that the supercritical branching random walk 
survives on large finite cubes with positive probability.
We adapt their result to the continuous-space case.

Our proof of \eqref{exp est}
relies on renormalization 
and  comparison with oriented percolation.
This scheme has been carried out for 
the contact process, see  e.g. 
Bezuidenhout and Grimmett \cite{BG90},
Durrett \cite{Dur91} or Liggett \cite{Lig99}.
Since in our case the geographic space is bounded
but the spin space is unbounded,
we use a different approach based on the genealogical structure.

The paper is organized as follows.
In Section 2 
we describe the model
and give our assumptions and results.
Sections 3 to 5 are devoted to proofs.

\section{The model, assumptions and results}

\emph{Description}.
The evolution of the system admits the following description. 
Each particle
``lives'' in $\Z ^\d$ (the discrete-space case)
or $\R ^\d$ (the continuous-space case) and
has two  exponential clocks with parameters $1$ and 
${\lambda}$, $\lambda > 1$.
When 
the first clock rings, the particle is deleted from the system (``death'').
When the second clock rings, the particle 
gives a birth to a new particle. 
After that the clocks are reset.
The offspring is distributed according 
to some radially symmetric dispersal kernel $a$.
Births outside of some  cube $\B$ are suppressed,
and there are no particles outside $\B$ at the beginning.

In the discrete-space case the state space of the process 
is $\Z _+ ^{\B}$, in the continuous-space case it is 
the collection of finite subsets of $\B$: $\{\eta \subset \B : |\eta \cap \B| < \infty \}$.
In either case we denote the state space by $\X$.

 The heuristic generator is given by 

  \[
  L F (\eta) = \sum\limits _{x \in \eta} \big\{ F(\eta \setminus \{ x\}) - F (\eta) \big\}
  + \lambda \sum\limits _{x \in \eta} \int\limits _{y \in X: x-y \in \B} 
  a(y-x) \big\{ F(\eta \cup \{ y\}) - F (\eta) \big\} \nu(dy),
 \]
where $\lambda > 0$ is the branching rate,
$F: \X \to \R _+$ is some function from an appropriate domain,
$X = \R ^\d$ and $\nu$ is the Lebesgue measure,
or $X = \Z ^\d$ and $\nu$ is the counting measure.
In both cases, 
\[
 \int\limits _{y \in X} a(y) \nu (dy) = 1.
\]

The process can be constructed in the following way.
Take a rooted tree $\E$ as in Figure 2.
To a vertex $\e$ we assign an independent vector
$(b_\e , d_\e , s _\e  )$ with values in 
$\R _+ \times \R _+ \times X$, where 
$X = \Z ^\d$ or $X = \R ^\d$. We take $b_\e$ and $d_\e$
to be exponentials with parameters $\lambda$ and $1$ respectively,
and $s _\e$ to be distributed according to $a$.
Assume that the particle to which $\e$
is assigned is born at time $t_\e$ at $x \in X$.
If $d_\e < b_\e$, the particle dies at time $t_\e + d_\e $,
otherwise the particle produces an offspring 
at time $t_\e + b_\e $. The position of the offspring is 
 $s _\e +x$. 
 The offspring is removed instantly
 if it is born outside $\B$.
The initial particle is assigned to the root of the tree.
This construction naturally allows us to endow 
the process with the genealogical structure.

\begin{figure}

  \begin{tikzpicture}
  
  \draw[thick] (-2,0) -- (0,0); 
  \draw[thick] (0,0) -- (4.5,2.5); 
  \draw[thick] (4.5,2.5) -- (10,3.5); 
  \draw[thick] (4.5,2.5) -- (8,1.5); 
   \draw[thick] (0,0) -- (1.5,-1);
  \draw[thick] (1.5,-1) -- (2.5,-0.333);
  \draw[thick] (1.5,-1) -- (5.5,-1.5);
  \draw[thick] (5.5,-1.5) -- (5.5,-2) node[anchor=north west] {$Q$};
  \draw[thick] (5.5,-1.5) -- (8.5,-1.5);
  \draw[thick] (8.5,-1.5) -- (10,-0.5);
  \draw[thick] (8.5,-1.5) -- (10,-2.5);
  
  \draw[thick,->] (-2,-3.5) -- (11,-3.5) node[anchor=north west] {time};

  \fill (2.5,-0.333) circle[radius=2pt];
  \fill (5.5,-2) circle[radius=2pt];
  \fill (8,1.5) circle[radius=2pt];
  
   \node[anchor=north ] at (10,-3.5) {$t_4$};   \draw[thick] (10,-3.4) -- (10,-3.6);
   \node[anchor=north ] at (0,-3.5) {$t_1$};   \draw[thick] (0,-3.4) -- (0,-3.6);
   \node[anchor=north ] at (-2,-3.5) {$0$};   
   \node[anchor=north ] at (2.5,-3.5) {$t_2$};   \draw[thick] (2.5,-3.4) -- (2.5,-3.6);
   \node[anchor=north ] at (5.5,-3.5) {$t_3$};   \draw[thick] (5.5,-3.4) -- (5.5,-3.6);

 \end{tikzpicture}
 \caption{  {\footnotesize Genealogical structure of the process. The first birth occurs at   \(t_1\), the first death at \(t_2\). The newly born at \(t_3\) particle is
 outside \(B\), so it dies instantly, hence the vertical line. There are \(3\) particles alive at \(t_4\).} } 
 \end{figure}
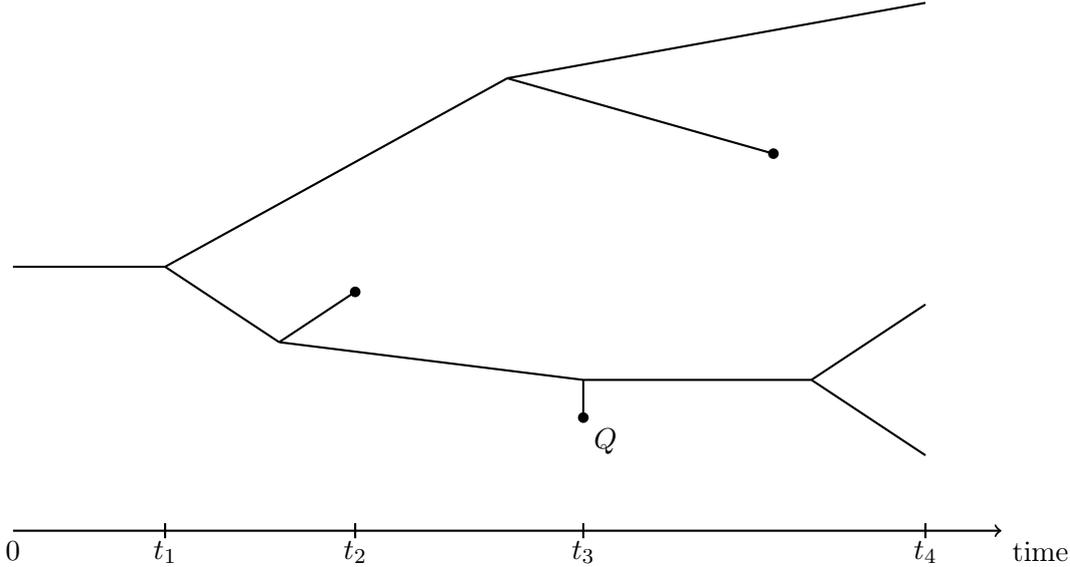

If $\beta$ is some collection of particles of the BRW alive at time $s$,
we denote by $(\eta ^{s, \beta} _{t})_{t \geq 0}$ the process 
starting from $\beta$ at $s$. Clearly, if $\alpha \subset \beta$,
then $\eta ^{s, \alpha} _{t} \subset \eta ^{s, \beta} _{t}$ for all $t \geq s$.
The process started from a single particle at $x \in \B$
is denoted by $(\eta _t ^{0,x})_{t \geq 0}$.

We write $(\eta _t ^{0,x})$ as a shorthand for $(\eta _t ^{0,x})_{t \geq 0}$,
meaning the whole trajectory of the process.
We say that the BRW \emph{survives on $\B$ with positive probability},
if there is an $x\in \B$ such that $P \{ (\eta _t ^{0,x}) \textrm{ survives} \} >0$.
Note that if the BRW survives on some cube with positive probability,
it also does so on a larger cube.

\emph{Assumptions and results}.
Let $a^{(n)}$ be the $n$-time convolution of $a$,
or the $n$-step transition function/density.
In the discrete-space case we say that $a$ is elliptic if (cf. \cite{BZ09})
 for any  $y \in \Z ^\d$
\begin{equation} \label{dally}
  a^{(n)}(y) > 0  \ \ \  \text{ for some } n \in \N.
\end{equation}
In the continuous-space case we say that $a$ is elliptic if 
for any $y \in \R ^\d$ and $r>0$,
\begin{equation} \label{solicit}
  \inf\limits _{z \in B(y,r)} a^{(n)}(z) > 0  \ \ \  \text{ for some } n \in \N,
\end{equation}
where $B(y,r)$ is the ball of radius $r$ around $y$.

We assume that $a$ is continuous (in discrete-space settings 
it amounts to no assumption)
and elliptic. 
Note that for the survival on a  cube we need
some kind of ellipticity of $a$:
for example, if $\d = 1$ and 
the support of $a$ is contained by $[1,\infty)$,
then the BRW dies out for every $\B$ and $\lambda > 0$.
 In the discrete-space settings,
the survival of the supercritical BRW ($\lambda >1$) 
on large cubes has been proven by
Mountford and Schinazi \cite{MS05}, 
for the BRW corresponding to the simple random walk,
and by
Bertacchi and Zucca \cite[Section 3]{BZ09}, under conditions 
similar to \eqref{dally} for a BRW on a general 
connected graph of bounded degree. 
The following theorem extends these results
to continuous-space settings.

\begin{thm}\label{cont space}
In the continuous-space case,
  the BRW survives on $\B$ with positive probability  provided
 that $\B$ is sufficiently large.
 
\end{thm}

  Let 
$\tau$ be the moment of extinction,
with convention that $\tau = \infty$ if the process survives.
Assume that $\B$ is sufficiently
large so that the process survives with positive probability.

For technical reasons, in the continuous-space case 
we will impose stronger conditions than \eqref{solicit}.
Let $\0$ be the origin in $\R _\d$,
$\Delta$ a 'cemetery' state, and
$ \tilde a _\B : \left( \B \cup \{\Delta \} \right) \times \mathscr{B}\left( \B \cup \{\Delta \} \right) 
 \to [0,\infty)$ be the transition function given by
 \begin{equation*}
  \tilde a _\B (x,B) = \int _{y \in B} a(y-x), \quad x, y \in \B,  B \in \mathscr{B}(\B),
 \end{equation*}
 $\tilde a _\B(x, \{ \Delta\}) = \int _{y \notin \B} a(y-x)$, and $\tilde a _\B (\Delta, \cdot) \equiv 0$.
 Here $\mathscr{B}(\B)$ is the collection of Borel subsets of $\B$.

First, assume that $P \{ (\eta _t ^{0,\0}) \textrm{ survives} \} >0$.
We further assume that
for every $r>0$ there exist
$N \in \N$ and $\tilde \delta >0$
such that
\begin{equation} \label{inane}
 \forall x \in \B \quad \sum\limits _{n=1} ^N  \tilde  a^{(n)} _\B (x,B(\0,r)) \geq \tilde \delta.
\end{equation}
and that
there is a small ball $B(\0,\bar r)$
such that for any $y \in B(\0,\bar r)$ and $\bar \delta >0$,
\begin{equation}\label{vitriol}
 P \{ (\eta _t ^{0,y}) \textrm{ survives} \} > 
 \bar \delta.
\end{equation}

Combining \eqref{inane} and \eqref{vitriol} gives the existence of $\delta > 0$ such that
\begin{equation} \label{counterfeit}
 \forall y \in \B \quad P \{ (\eta _t ^{0,y}) \textrm{ survives} \} >\delta.
\end{equation}

The following theorem is the main result of this paper.

\begin{thm} \label{main thm}
Under the above assumptions,
 \eqref{exp est} holds.
\end{thm}

\begin{rmk}\label{resuscitate}
 Assumption \eqref{vitriol} is not very restrictive due to 
 the following observation. Assume that 
$P \{ (\eta _t ^{0,\0}) \textrm{ survives} \} = p_{_\B} >0$
and let $l$ be the length of an edge of $\B$.
Then for a cube  $ \B ^ \varepsilon$ with the edge length 
 $l + 2\varepsilon$, $\varepsilon >0$,
 and  for all $y \in (- \varepsilon, \varepsilon) ^\d$
\begin{equation*}
 P \{ (\eta _t ^{0,y}) \textrm{ survives} \} \geq
 p_{_\B}.
\end{equation*}

\end{rmk}

\begin{rmk}
For the supercritical 
process on the whole space, $\Z ^\d$ or $\R ^\d$,
\eqref{exp est}
comes down to the corresponding estimate
for the Galton--Watson process,
since $X _t := |\eta _t|$ is a birth-death process with transition rates
\begin{gather*}
 n \to n+1  \ \ \ \textrm{ at rate } \lambda n,
 \\
 n \to n-1  \ \ \ \textrm{ at rate }  n.
\end{gather*}
\end{rmk}

\section{Proof of Theorem \ref{cont space}}

The idea is to couple a continuous-space supercritical BRW
with a discrete-space one and then use the result of 
 \cite{BZ09}. With no loss of generality
we assume that the length of an edge of $\B$ is 
a natural number.

For $n \in \N$ and 
$j = (j_1,...,j_\d) \in \frac{1}{2^n} \Z ^\d \cap \textrm{int}(\B) $,
where $\textrm{int}(\B)$ is the interior of $\B$,
we define 
\[
a_n(j) = \frac{1}{2^{n\d}} \inf \{ a(x-y): x \in [-\frac{1}{2^{n+1}},\frac{1}{2^{n+1}})^{\d}, 
  y \in \prod\limits _{k=1} ^\d [j_k-\frac{1}{2^{n+1}},j_k+\frac{1}{2^{n+1}}) \}.
\]
Note that $a_n$ is elliptic.
Since $a$ is continuous, we have 
\begin{equation}
 \sum\limits _{j \in \frac{1}{2^n} \Z ^\d \cap \textrm{int}(\B)}
  a_n (j)  \to \int\limits _{\R ^\d} a(x) dx,
\end{equation}
therefore $\sum\limits _{j \in \frac{1}{2^n} \Z ^\d}
  a_n (j) > 1$ for sufficiently large $n$.
  We will choose such an $n \in \N$ 
  and couple the given continuous-space BRW $(\eta _t)$
  with
  discrete-space BRW $(\eta ^{(n)} _t)$ 
  on $\frac{1}{2^n} \Z ^\d$ with kernel $a_n$ as follows.
  Each particle $q$ from $(\eta ^{(n)} _t)$  is associated
  to a particle $s(q)$ from $(\eta  _t)$, and 
  no particle from $(\eta  _t)$ may have two 
  particles from $(\eta ^{(n)}  _t)$ associated to it,
  so that $s :\eta ^{(n)} _t \to  \eta  _t $ is an injection for each $t$.
  We consider $(\eta _t)$ started from one particle at the origin.
  We let
  $\eta ^{(n)} _0$ to have one particle at the origin of 
  $\frac{1}{2^n} \Z ^\d$, which we associate to 
  the initial particle of $(\eta _t)$.

  If a particle $s(q)$ at $x$ gives birth to a new particle at $y$ at a time $s$, 
  where $x \in [j^x_k-\frac{1}{2^{n+1}},j^x_k+\frac{1}{2^{n+1}})$, 
  $y \in [j^y_k-\frac{1}{2^{n+1}},j^y_k+\frac{1}{2^{n+1}})$ for 
  some $j^x, j^y \in \frac{1}{2^n} \Z ^\d$,
 then the associated to the parent particle $q$ at $j^x$ gives birth
 to a new particle at $j^y$ with probability $\frac{a_n(j^y - j^x)}{a(y-x)}$,
 provided that the  particle  $s(q)$ exists and is alive.
  We associate the  newborn particles to each other. Also, 
 associated particles die simultaneously.
  
 It is clear that $|\eta ^{(n)} _t| \leq |\eta  _t|$ for all $t \geq 0$;
 in particular, if $(\eta ^{(n)} _t)$ survives, then so does $(\eta  _t)$.
 It remains to note that
  from  \cite[Theorem 3.1]{BZ09} we know that 
 $(\eta ^{(n)} _t)$ survives on a sufficiently large 
 finite cube with positive probability.
 \qed

\section{Proof of Theorem \ref{main thm}}

We prove Theorem \ref{main thm} concurrently in
discrete and continuous settings,
because the ideas involved 
are very similar.
We endow our system with the
genealogical structure, so that we can talk about 
ancestors and descendants.
Without loss of generality we assume
that $\B$ is centered at the origin.
Furthermore, we assume without loss of generality
that in the discrete-space case
the random walk on $\B$ with the kernel $a_\B$
is irreducible. Here for $x,y \in \B$
\[
 a_\B (y,x) = a(y-x) + I\{ x = y \} \sum\limits _{z \ne \B} a (z-x).
\]
Concerning the last assumption, see Remark \ref{innuendo}.

\begin{lem} \label{trite}
In the discrete-space case,
 for any $\varepsilon > 0$ there are 
 $T>0$ and $M \in \N$ such that 
 \begin{equation}
 c_{ij} := P \{ \eta ^{0, M I _{A_i}} _T \geq M I _{A_j} \} \geq 1 - \varepsilon, \quad i,j = 1,2,
 \end{equation}
 where 
 $A_1 = \{(x_1,...,x_\d) \in \B \mid x_1 \geq 0 \}$ and $A _2 = \B \setminus A_1$.
 
\end{lem}

\textbf{Proof}. The BRW can be considered as a 
continuous-time Markov chain on $\B ^{\Z _+}$.
Since zero state is a trap that can be reached from any state, 
any finite subset of $\B ^{\Z _+}$ is transient. In particular,
for any
$L>0$
\begin{equation} \label{pull in horns}
   P \{ \tau = \infty, \max\limits _{x \in \B} \eta _t (x) \leq L \} \to 0,  \ \ \ t \to \infty.
\end{equation}
Let us choose $M$ so large that 
\[
 P \{  \eta ^{0, M I_{A_i}}  \text{ dies out} \} \leq 1 - \frac{\varepsilon}{4}, \ \ \ i = 1,2.
\]

Proceeding further,
let
us choose $L$ so large that the following is satisfied:
for any $x \in \B$, process started at $0$ 
from $L$ particles in $x$ has at time $1$ at least $M$ particles
everywhere on $\B$ with probability larger than $1 - \frac{\varepsilon}{4}$.
Choosing now $T $ so large that 
\[
 P \{ \max\limits _{x \in \B} \eta ^{0, M I_{A_i}} _{T-1} (x) \geq L \} \geq 1 - \frac{\varepsilon}{2}, \ \ \ i = 1,2.
\]
completes the proof.
\qed

\begin{rmk}\label{innuendo}
It can be that the random walk with transition function $a_\B$
is not irreducible on $\B$. As an example,
let us take $\d = 1$, $\B = \{-2,...,2\}$ and  $a(x) = \frac 12 I\{ |x| = 5\} $
and note that the corresponding BRW survives with positive probability if $\lambda >2$.
If this is the case,
there is a component $\bar \B \subset \B$
such that the BRW started from a single particle 
in $\bar \B \subset \B$ survives with positive probability
within $\bar \B$ (that is, with births outside $\bar \B$
being suppressed; in the above example
$\bar \B$ would be $\{ -2, 2 \}$). 
The above  lemma still holds
provided that $A _i$ is replaced by 
$ A _i \cap \bar \B$, $i = 1,2$.
\end{rmk}

Define
$$
Q _+ = \B \cap \left\{ \x \in \R ^\d: \x = (x_1,...,x_\d)
\text{ with } x_1 \geq 0   \right\}$$
and 
$$
Q _- = \B \cap \left\{ \x \in \R ^\d: \x = (x_1,...,x_\d)
\text{ with } x_1 < 0  \right\}.
$$
For $M \in \N$,
let
\[
 A ^M _+ = \{ \eta \in \Gamma _0 (\B): |\eta \cap Q_+| >M  \}
\]
and \[
 A ^M _- = \{ \eta \in \Gamma _0 (\B): |\eta \cap Q_-| >M  \}.
\]

\begin{lem} \label{glitch}
 In the continuous-space case, 
 for any $\varepsilon > 0$
 there are 
 $T>0$ and $M \in \N$ such that
  \begin{equation}
  P \{ \eta ^{0, \eta _0} _T \in A ^M _j \} \geq 1 - \varepsilon
 \end{equation}
 for 
 any $\eta _0 \in A ^M _i$. Here each of the indices $i$ and $j$ can be either $+$ or $-$.

\end{lem}

\textbf{Proof}. By a similar argument,
for any $n \in \N$ the set $\{\eta \subset \B : |\eta| = n \}$
is transient
in the sense that a.s. 
it is entered finitely many times only.
The counterpart of \eqref{pull in horns} is 
\begin{equation*}
   P \{ \tau = \infty, | \eta _t (x)| \leq L \} \to 0,  \ \ \ t \to \infty.
\end{equation*}
By \eqref{counterfeit},
the probability of survival is separated from $0$.
We can choose $M$ so large that 
\[
 P \{  (\eta ^{0, \eta _0} _t)  \text{ dies out} \} \leq 1 - \frac{\varepsilon}{4}
\]
for any $\eta _0 \in A ^M _i$, $i = +,-$,
then $L$ so large that any
process started from $L$ particles at time $0$
is in the intersection $A ^M _+ \cap A ^M _-$
by time $1$ with high probability ($1 - \frac{\varepsilon}{4}$ is sufficient),
and finally we choose $T$ so that
\[
 P \{  |\eta ^{0, \eta _0} _{T-1}| \geq L   \} \leq 1 - \frac{\varepsilon}{2}
\]
for any $\eta _0 \in A ^M _i$, $i = +,-$,
and the proof goes as in Lemma \ref{trite}.
\qed

Let $G = \{ (n,m) : n+m \textrm{ is even} \}$. We will
use Lemmas \ref{trite} and \ref{glitch}
to make a comparison with the
 oriented percolation process 
on $G$. Let $(n,m)$ be connected to 
$(n+1, m+1)$ and $(n-1, m+1)$. 
Each bond is open with 
probability $p$ independently 
of the other bonds. 
We say that percolation occurs
if there is an infinite path starting from 
the origin.
The model is well-known, 
see e.g. Durrett \cite{Dur84, Dur88}.

Let $p_c$ be the critical value for 
independent oriented percolation in two dimension,
and let 
$$
\sigma = \min \big\{ m \in \N : \textrm{ there is no open path from } (0,0) \textrm{ to } \{(k,m)\mid k \in \Z \} \big\},
$$
the moment of extinction of the percolation process.
We use the following estimate in the proof of Theorem \ref{main thm}.

\begin{lem}[\protect{\cite{Dur84}}]
Assume that $p > p_c$.
 Then there are $q_1, C_1 >0$ such that
 \begin{equation} \label{veracity}
  P \{ r < \sigma < \infty  \} \leq C_1 e^{-q_1 r},  \quad  r \geq 0.
   \end{equation}

\end{lem}

\textbf{Proof of Theorem} \ref{main thm} \textbf{in the discrete-space case}.
Let us take $M$ and $T$ so large that 
Lemma \ref{trite} is satisfied 
with $1 - \varepsilon = p > p_c$.
Let $(u_n)_{n\in G}$
be a sequence of independent random variables
distributed uniformly on $[0,1]$,
independent of everything introduced so far.
Denote also 
\[
 c_{ij} = P \{ \eta ^{0, M I _{A_i}} _T \geq M I _{A_j}  \} \geq p.
\]

Let $\tau _1 = \tau \wedge  \inf \{t : \eta _t \geq M I _{A_2} \}$.
Since every particle alive at some time $t_0$ 
produces
by the time $t_0 + 1$ 
so many particles as to dominate $M I _{A_1}$
with positive probability
separated from zero, 
$\tau _1$ is dominated by a geometric random variable
and has subexponential tails (see \mbox{Section  \ref{subexp tails}}
for the precise meaning of ``subexponential tails'').  
If the process does not die out at $\tau _1$, 
then we build an oriented bond percolation process on $G$
 according to the following procedure.

Choose  a collection of particles $\alpha _{(0,0)}$
alive at time $\tau _1$
in such a way that $S(\alpha_{(0,0)}) = M I _{A_2}$.
Here $S(\alpha_{(0,0)}) = M I _{A_2}$ means
that $\alpha_{(0,0)}$ has exactly $M$
particles at every site from $A_2$
and has no particles outside $A_2$.
In our construction, $S(\alpha_{(n,m)}) = M I _{A_2}$
if $m \equiv n \mod 4 $, and $S(\alpha_{(n,m)}) = M I _{A_1}$
if $m \equiv n + 2 \mod 4 $.

We say the edge $\langle (0,0), (1,1) \rangle$ from $(0,0)$ to $(1,1)$ is open
if both
$$
\{ \eta ^{\tau _1, \alpha_{(0,0)}} _{\tau _1 +T} \geq M I _{A_2} \}
$$
and 
$$
\{u_{\langle (0,0), (1,1) \rangle} < \frac{p}{c_{22}}\}
$$
occur, and we say that the edge $ \langle (0,0), (-1,1) \rangle $ is open if both
$\{ \eta ^{\tau _1, \alpha_{(0,0)}} _{\tau _1 +T} \geq M I _{A_1} \}$
and $\{u_{\langle (0,0), (-1,1) \rangle} < \frac{p}{c_{21}}\}$
occur. If $\langle (0,0), (1,1) \rangle$ is open, then we choose 
$\alpha _{(1,1)}$ in such a way that 
$S(\alpha_{(1,1)}) = M I _{A_2}$
and that every particle from $\alpha _{(1,1)}$
is an descendant of a particle from $\alpha _{(0,0)}$
(here we consider a particle to be a descendant of itself
provided that it is still alive).
Similarly, if $\langle (0,0), (-1,1) \rangle$ is open, we choose 
$\alpha _{(-1,1)}$ in such a way that 
$S(\alpha_{(-1,1)}) = M I _{A_1}$
and that every particle from $\alpha _{(-1,1)}$
is an descendant of a particle from $\alpha _{(0,0)}$.
Further proceeding, assume that there is an open path from
the origin to $(n,m)$,
and a collection $\alpha_{(n,m)}$ of particles 
alive at $\tau _1 + mT$ is chosen, such that
\begin{equation}\label{berate}
  S(\alpha_{(n,m)}) =
  \begin{cases} 
      \hfill M I _{A_1}    \hfill & \text{ if $m \equiv n+2 \mod 4 $,} \\
      \hfill M I _{A_2}    \hfill & \text{ if $m \equiv n \mod 4 $}. \\
  \end{cases}
\end{equation}

For $m \equiv n \mod 4 $, 
we let $\langle (n,m),(n+1,m+1) \rangle$ be open if 
$\{ \eta ^{\tau _1 + mT, \alpha_{(n,m)}} _{\tau _1 + (m+1)T} \geq M I _{A_2} \}$
and $\{u_{\langle (n,m), (n+1,m+1) \rangle} < \frac{p}{c_{22}}\}$
occur, and $\langle (n,m),(n-1,m+1) \rangle$ is open if 
$\{ \eta ^{\tau _1 + mT, \alpha_{(n,m)}} _{\tau _1 + (m+1)T} \geq M I _{A_1} \}$
and $\{u_{\langle (n,m), (n-1,m+1) \rangle} < \frac{p}{c_{21}}\}$ do.
Similarly, for $m \equiv n+2 \mod 4 $,
$\langle (n,m),(n+1,m+1) \rangle$ is open 
if $\{ \eta ^{\tau _1 + mT, \alpha_{(n,m)}} _{\tau _1 + (m+1)T} \geq M I _{A_1} \}$
and $\{ u_{\langle (n,m), (n+1,m+1) \rangle} < \frac{p}{c_{11}} \}$
occur, 
and $\langle (n,m),(n-1,m+1) \rangle$ is open if 
$\{ \eta ^{\tau _1 + mT, \alpha_{(n,m)}} _{\tau _1 + (m+1)T} \geq M I _{A_2} \}$
and $\{ u_{\langle (n,m), (n-1,m+1) \rangle} < \frac{p}{c_{12}} \} $ do.
Furthermore, if $\langle (n,m),(n \pm 1,m+1) \rangle$ is open,
we choose $\alpha_{(n \pm 1,m+1)}$ in such a way that each particle
from $\alpha_{(n \pm 1,m+1)}$ is a descendant from a particle
from $\alpha_{(n,m)}$ 
and \eqref{berate} is satisfied.

If there is no open path to $(n,m)$,
then $\alpha_{(n,m)}$ is not defined, and
we may take $\langle (n,m),(n \pm 1, m+1) \rangle$
to be open iff $u_{\langle (n,m), (n \pm 1,m+1) \rangle} < p$.

Thus we get the desired percolation process,
in which edges are open independently with 
probability $p$,
and
which
is constructed 
in such a way that percolation implies survival of $(\eta _t) _{t \geq 0}$.
Let $\sigma _1$ be the lifetime 
of the percolation process.
If percolation doesn't occur
but the BRW  still lives,
we start anew and 
on $\{\tau > \tau _1, \sigma _1 < \infty \}$
define $\tau _2$
analogously to $\tau _1$, 
$$
\tau _2 = \tau \wedge \inf \{t > \tau _1 +  \sigma _1 T: \eta _t \geq M I_{A_2} \}.
$$ 
If, after some time, the BRW dies out at some $\tau _i$,
then we use an independent collection 
of oriented percolation processes
to define $\sigma _i$ until 
the first time percolation occurs.
Let $g \in \N $ the 
number of the first percolation process that survives, that is $\sigma _{g-1} < \infty$ 
and $\sigma _{g-1} = \infty$. Clearly, $g$ has a geometric distribution.
A.s. on $\{\tau < \infty \}$ we have
$$
\tau \leq I\{ g \ne 1 \} \sum\limits _{j=1 } ^{g-1} (\tau _j + \sigma _j T ) + \tau _g, 
$$
 where $\tau _j$,  $\sigma _j$ have subexponential tails
and $g$ has a geometric distribution. 
It remains to apply two lemmas from Section \ref{subexp tails}.
\qed

\textbf{Proof of Theorem} \ref{main thm} \textbf{in the continuous-space case}.
We will use a similar percolation argument 
to prove Theorem \ref{main thm} in continuous-space settings.
Take $T >0$ and $M \in \N$ so large that
Lemma \ref{glitch} is satisfied with
$1 - \varepsilon = p \in (p_c,1)$. 
Similarly to the discrete-space case,
let $\tau _1 = \tau \wedge  \inf \{t : \eta _t \in A ^{M}_- \}$.
If $\tau _1 \ne \tau$,
choose a minimal $\alpha_{(0,0)}$
such that $\alpha_{(0,0)} \subset \eta _{\tau _1}$ 
and $\alpha_{(0,0)} \in A ^{M}_-$.

Let $\bar \alpha_{(0,0)}$ be some collection of particles
alive at time $0$ and having spatial positions
identical to particles from $ \alpha_{(0,0)}$.
We declare $\langle (0,0), (-1,1) \rangle$ to be open
if
$\{ \eta ^{\tau _1, \alpha_{(0,0)}} _{\tau _1 +T} \in A ^{M}_- \}$
and
$$u_{\langle (0,0), (-1,1) \rangle} < p 
\big(
P 
\{ \eta ^{0,\bar \alpha_{(0,0)} } _{ T} \in A ^{M}_- \} \big) ^{-1},$$
and so on, proceeding exactly as in the discrete-space case.
That will yield the desired result.
\qed
\begin{rmk}
  In the proof of Theorem \ref{main thm}
we tacitly 
assumed that the
  strong Markov property holds at $\tau _1, \tau _2, ...$.
We could prove that $(\eta _t)$ has the strong Markov property,
but in this case it is easier to replace $\tau _1 $ with 
\[
 \tilde \tau _1 = \left \lceil{\tau}\right \rceil \wedge 
\min \{n \in \N : \eta _n \geq M I _{A_2} \},
 \]
where $ \left \lceil{\cdot}\right \rceil$
is the ceiling function, and use the fact that 
the strong Markov property is satisfied 
for the stopping times
which take countably many values only,
see e.g. Kallenberg \cite[Proposition 8.9]{KallenbergFound}.
In a similar way we can replace $\sigma _1 , \tau _2$, and so on.
The proof needs no further changes.
\end{rmk}

\section{Subexponential tails}\label{subexp tails}

We say that a random variable $X$ has subexponential tails if there are 
$C_{_X}, q_{_X} > 0 $ such that 
\[
 P \{X \geq x \} \leq C_{_X} e ^{-q_{_X} x}, \ \ \ x \geq 0.
\]
Note that $E e^{\theta X} < \infty$ if $\theta < q_{_X}$.

\begin{lem}
 Let $X$ and $Y$ be independent random variables with 
 subexponential tails. Then their sum has subexponential tails too.
\end{lem}
\textbf{Proof}.
$P \{X + Y \geq 2 z \} \leq P \{X \geq  z \} + P \{ Y \geq z \}$.
\qed
\begin{lem}
 Let $X_1, X_2, ...$ be a sequence of i.i.d. random variables with subexponential tails,
 and let $g$ be an independent random variable with a geometric distribution,
 \[
  P \{ g = m \} = (1- p) p^{m - 1}, \ \ \ m \in \N,
 \]
 where $p \in (0,1)$.
  Then $S = \sum\limits _{j=1} ^g X _j$ has subexponential tails.
\end{lem}
\textbf{Proof}.
 By the Lebesgue dominated convergence theorem  there exists $\theta >0 $ such that $Ee^{\theta X_1} < \frac 1p$.
 For such  a  $\theta$,
 \[
  E e ^{\theta S} = \sum\limits _{m =1 } P \{g = m \} (Ee^{\theta X_1}) ^{m} < \infty,
 \]
hence by Chebyshev's inequality
\[
 P \{ S > x  \} \leq {E e ^{\theta S}}e^{-\theta x}.
\]

 \section*{Acknowledgement}

 I would like to thank Yuri Kondratiev 
 and Tyll Kr\" uger for the
 inspiring discussions. 
This paper was partially supported by the DFG through the SFB 701
``Spektrale Strukturen und Topologische Methoden in der Mathematik'' and
by the European Commission under the project STREVCOMS PIRSES-2013-612669.

\bibliographystyle{alpha}
\bibliography{Sinus}

\begin{thebibliography}{Dur91}

\bibitem[BG90]{BG90}
C.~Bezuidenhout and G.~Grimmett.
\newblock The critical contact process dies out.
\newblock {\em Ann. Probab.}, 18(4):1462--1482, 1990.

\bibitem[BZ09]{BZ09}
D.~Bertacchi and F.~Zucca.
\newblock Approximating critical parameters of branching random walks.
\newblock {\em J. Appl. Probab.}, 46(2):463--478, 2009.

\bibitem[BZ15]{BZ15}
D.~Bertacchi and F.~Zucca.
\newblock Branching random walks and multi-type contact-processes on the
  percolation cluster of {$\mathbb{Z}^d$}.
\newblock {\em Ann. Appl. Probab.}, 25(4):1993--2012, 2015.

\bibitem[Dur84]{Dur84}
R.~Durrett.
\newblock Oriented percolation in two dimensions.
\newblock {\em Ann. Probab.}, 12(4):999--1040, 1984.

\bibitem[Dur88]{Dur88}
R.~Durrett.
\newblock {\em Lecture notes on particle systems and percolation}.
\newblock Wadsworth \& Brooks/Cole Advanced Books \& Software, Pacific Grove,
  CA, 1988.

\bibitem[Dur91]{Dur91}
R.~Durrett.
\newblock The contact process, 1974--1989.
\newblock In {\em Mathematics of random media}, volume~27 of {\em Lectures in
  Appl. Math.}, pages 1--18. Amer. Math. Soc., Providence, RI, 1991.

\bibitem[Kal02]{KallenbergFound}
O.~Kallenberg.
\newblock {\em Foundations of modern probability}.
\newblock Probability and its Applications. Springer-Verlag, second edition,
  2002.

\bibitem[Lig99]{Lig99}
T.~M. Liggett.
\newblock {\em Stochastic interacting systems: contact, voter and exclusion
  processes}, volume 324.
\newblock Springer-Verlag, 1999.

\bibitem[MS05]{MS05}
T.~Mountford and R.~B. Schinazi.
\newblock A note on branching random walks on finite sets.
\newblock {\em J. Appl. Probab.}, 42(1):287--294, 2005.

\end{thebibliography}

\end{document}